\documentclass[adm]{birkmult}
\usepackage{amsmath}
\usepackage{amsfonts}
\usepackage{amssymb}
\providecommand{\U}[1]{\protect\rule{.1in}{.1in}}

\newtheorem{theorem}{Theorem}[section]

\newtheorem{lemma}[theorem]{Lemma}
\newtheorem{proposition}[theorem]{Proposition}
\theoremstyle{definition}

\newtheorem*{acknowledgment}{Acknowledgment}

\newtheorem{note}[theorem]{Note}
\numberwithin{equation}{section}
\begin{document}
\author[D.I. Dais]{Dimitrios I. Dais}
\address{University of Crete, Department of Mathematics, Division Algebra and Geometry,
Knossos Avenue, P.O. Box 2208, GR-71409, Heraklion, Crete, Greece}
\email{ddais@math.uoc.gr}
\author[B. Nill]{Benjamin Nill}
\address{Freie Universit\"{a}t Berlin, Institut f\"{u}r Mathematik, Arbeitsgruppe
Gitterpolytope, Arnim\-allee 3, 14195 Berlin, Germany}
\email{nill@math.fu-berlin.de}
\date{}
\title[Boundedness for Toric Log Del Pezzo Surfaces]{A Boundedness Result \\for Toric Log Del Pezzo Surfaces}

\begin{abstract}
In this paper we give an upper bound for the Picard number of the rational
surfaces which resolve minimally the singularities of toric log Del Pezzo
surfaces of given index $\ell$. This upper bound turns out to be a quadratic
polynomial in the variable $\ell$.

\end{abstract}
\maketitle

\section{Introduction\label{INTRO}}

\noindent{}A normal complex surface $X$ with at worst log terminal
singularities, i.e., quotient singularities, is called \textit{log Del Pezzo}
\textit{surface} if its anticanonical divisor $-K_{X}$ is a $\mathbb{Q}%
$-Cartier ample divisor. The \textit{index} of such a surface is defined to
be the smallest positive integer $\ell$ for which $\ell K_{X}$ is a Cartier
divisor. Every log Del \ Pezzo surface is isomorphic to the
\textit{anticanonical model} (in the sense of Sakai \cite{Sakai}) of the
rational surface obtained by its minimal desingularization. The following
Theorem is due to Nikulin \cite{Nikulin2} (for related results cf. \cite{Alex, Nikulin3}):

\begin{theorem}
\label{NIKTHM}Let $X$ be a log Del Pezzo surface of index $\ell$ and
$\widetilde{X}\longrightarrow X$ be its minimal desingularization. Then the
Picard number $\rho(\widetilde{X})$ of $\widetilde{X}$ \emph{(}i.e., the rank
of its Picard group\emph{)} is bounded by%
\begin{equation}
\rho(\widetilde{X})<c\cdot\ell^{\frac{7}{2}}, \label{NIKULINSUB}%
\end{equation}
where $c$ is an absolute constant.
\end{theorem}

\noindent{}The \textit{toric} log Del Pezzo surfaces, i.e., those which are
equipped with an algebraic action of a $2$-dimensional algebraic torus
$\mathbb{T},$ and contain an open dense $\mathbb{T}$-orbit, constitute a
special subclass within the entire class of all log Del Pezzo
surfaces.\ \ (For instance, in the toric case, only \textit{cyclic} quotient
singularities can occur.) To indicate how these two classes differ in
practice, it would be enough to recall some known results for log Del Pezzo
surfaces with Picard number $=1$ and index $\ell\leq2$:\smallskip

\noindent{}(i) Excluding the \textquotedblleft exceptional\textquotedblright%
\ $2D_{4}$-case, there exist, up to isomorphism, exactly $30$ surfaces of this
kind having index $\ell=1$ (see \cite[Thm. 4.3]{Alex-Nik} or \cite[Thm.
1.2]{Ye}). Among them there are $16$ having at worst cyclic quotient
singularities. By \cite[Thm. 6.10]{Dais} we see that only $5$ out of these
$16$ surfaces are toric (associated to the $5$ \textit{reflexive}
triangles).\smallskip

\noindent{}(ii) Up to isomorphism, there exist exactly $18$ surfaces of this
kind having index $\ell=2$ (see \cite[Thm. 4.2]{Alex-Nik} or \cite[Thm. 1.1
(1)]{Kojima}). Among them there are $14$ having only cyclic quotient
singularities. By \cite[Thm. 6.12]{Dais} we see that only $7$ out of these
$14$ surfaces are toric.

The purpose of this paper is to prove an analogue of (\ref{NIKULINSUB}) for
toric log Del Pezzo surfaces of given index.

\begin{theorem}
\label{main}Let $X_{Q}$ be a toric log Del Pezzo surface of index $\ell$
\emph{(}associated to the lattice polygon $Q$\emph{) }and $\widetilde{X}%
_{Q}\longrightarrow X_{Q}$ be its minimal desingularization. Then
$\rho(\widetilde{X}_{Q})$ is bounded as follows\emph{:}%
\begin{equation}
\rho(\widetilde{X}_{Q})\leq\left\{
\begin{array}
[c]{ll}%
7, & \text{\emph{if} }\ell=1,\\
8 \ell^{2} - 6 \ell+ 3, & \text{\emph{if} }\ell\geq2.
\end{array}
\right.  \label{mainbound}%
\end{equation}

\end{theorem}

\noindent{}Our proof uses tools from toric and discrete geometry.

\begin{acknowledgment}{\rm The second author is a member of the Research
Group Lattice Polytopes, led by Christian Haase and supported by Emmy Noether fellowship HA 4383/1
of the German Research Foundation (DFG).}
\end{acknowledgment}

\section{Toric log Del Pezzo surfaces\label{PRELIM}}

\noindent{}Let $Q\subset\mathbb{R}^{2}$ be a (convex) polygon. Denote by
$\mathcal{V}(Q)$ and $\mathcal{F}(Q)$ the set of its vertices and the set of
its facets (edges), respectively. $Q$ will be called an \textit{LDP-polygon}
if it contains the origin in its interior, and its vertices belong to
$\mathbb{Z}^{2}$ and are primitive. If $Q$ is an LDP-polygon, we shall denote
by $X_{Q}$ the compact toric surface constructed by means of the fan
\[
\Delta_{Q}:=\left\{  \left.  \text{the cones }\sigma_{F}\, \text{\ together
with their faces\ }\right\vert \,F\in\mathcal{F}(Q)\right\}  ,
\]
where $\sigma_{F}:=\left\{  \left.  \lambda\mathbf{x}\, \right\vert \,
\mathbf{x}\in F\text{ and }\lambda\in\mathbb{R}_{\geq0}\right\}  $ for all
$F\in\mathcal{F}(Q).$ It is known (cf. \cite[Remark 6.7]{Dais}) that every
toric log Del Pezzo surface is isomorphic to an $X_{Q},$ for a suitable
LDP-polygon $Q.$ Moreover, every cone $\sigma_{F}$ is lattice-equivalent to
the cone $\mathbb{R}_{\geq0}\binom{1}{0}+\mathbb{R}_{\geq0}\binom{p_{F}}%
{q_{F}},$ for suitable relatively prime integers $p_{F},q_{F},$ with $0\leq
p_{F}<q_{F}.$ (These are uniquely determined, up to replacement of $p_{F}$ by
its \textit{socius} $\widehat{p}_{F},$ i.e., by the integer $\widehat{p}_{F},$
$0\leq\widehat{p}_{F}<q_{F},$ satisfying gcd$(\widehat{p}_{F},q_{F})=1$ and
$p_{F}\widehat{p}_{F}\equiv1($mod $q_{F}).$) The affine toric variety
$U_{F}:=$ Spec$\left(  \mathbb{C}[\sigma_{F}^{\vee}\cap(\mathbb{Z}^{2})^{\vee
}]\right)  $ (where $\sigma_{F}^{\vee}$ denotes the dual cone of $\sigma_{F}$
and $(\mathbb{Z}^{2})^{\vee}$ the dual lattice of $\mathbb{Z}^{2}$) is
$\cong\mathbb{C}^{2}$ only if $q_{F}=1.$ Otherwise, the orbit orb$(\sigma
_{F})\in U_{F}$ of $\sigma_{F},$ i.e., the single point remaining fixed under
the canonical action of the algebraic torus $\mathbb{T}:=$ Hom$_{\mathbb{Z}%
}((\mathbb{Z}^{2})^{\vee},\mathbb{C}^{\star})$ on $U_{F}$, is a cyclic
quotient singularity. In particular, $U_{F}\cong\mathbb{C}^{2}/G_{F}=$
Spec$(\mathbb{C}[z_{1},z_{2}]^{G_{F}}),$ with $G_{F}\subset$ GL$\left(
2,\mathbb{C}\right)  $ denoting the cyclic group of order $q_{F}$ which is
generated by diag$(\zeta_{q_{F}}^{-p_{F}},\zeta_{q_{F}})$ (for $\zeta_{q_{F}}$
a $q_{F}$-th root of unity). Hence, the singular locus of $X_{Q}$ equals
\[
\text{Sing}(X_{Q})=\left\{  \left.  \text{orb}(\sigma_{F})\right\vert \,F\in
I_{Q}\right\}  ,
\]
where $I_{Q}:=\left\{  \left.  F\in\mathcal{F}(Q)\right\vert \,q_{F}%
>1\right\}  .$ Its subset $\{ \left.  \text{orb}(\sigma_{F})\right\vert
\,F\in\breve{I}_{Q}\},$ with $\breve{I}_{Q}$ defined to be $\breve{I}%
_{Q}:=\left\{  \left.  F\in I_{Q}\, \right\vert \,p_{F}=1\right\}  ,$ is the
set of the \textit{Gorenstein singularities} of $X_{Q}.$

The minimal desingularization of the surface $X_{Q}$ can be described as
follows: Equip the minimal generators of $\Delta_{Q}$ with an order (e.g.,
anticlockwise), and assume that for every $F\in\mathcal{F}(Q)$ the cone
$\sigma_{F}$ has $\mathbf{n}^{(F)},\mathbf{n}^{\prime(F)}\in\mathbb{Z}^{2}$ as
minimal generators ($\sigma_{F}=\mathbb{R}_{\geq0}\, \mathbf{n}^{(F)}%
+\mathbb{R}_{\geq0}\, \mathbf{n}^{\prime(F)}$), with $\mathbf{n}^{(F)}$ coming
first w.r.t. this order. Next, for all $F\in I_{Q},$ consider the
negative-regular continued fraction expansion of
\begin{equation}
\frac{q_{F}}{q_{F}-p_{F}}=\left[  \! \! \left[  b_{1}^{(F)},b_{2}^{(F)}%
,\ldots,b_{s_{F}}^{(F)}\right]  \! \! \right]  :=b_{1}^{(F)}%
-\cfrac{1}{b_{2}^{(F)}-\cfrac{1}{\begin{array} [c]{cc}\ddots & \\ & -\cfrac{1}{b_{s_{F}}^{(F)}}\end{array} }}\ \ ,
\label{EXPPQCF}%
\end{equation}
and define $\mathbf{u}_{0}^{(F)}:=\mathbf{n}^{(F)},$ $\mathbf{u}_{1}%
^{(F)}:=\frac{1}{q_{F}}((q_{F}-p_{F})\mathbf{n}^{(F)}+\mathbf{n}^{\prime
(F)}),$ and lattice points $\{ \mathbf{u}_{j}^{(F)}\left\vert \,2\leq j\leq
s_{F}+1\right.  \}$ by the formulae
\[
\mathbf{u}_{j+1}^{(F)}:=b_{j}^{(F)}\mathbf{u}_{j}^{(F)}-\mathbf{u}_{j-1}%
^{(F)},\ \ \forall j\in\{1,\ldots,s_{F}\}.\
\]
It is easy to see that $\mathbf{u}_{s_{F}+1}^{(F)}=\mathbf{n}^{\prime(F)},$
and that the integers $b_{j}^{(F)}$ are $\geq2,$ for all $j\in\{1,\ldots
,s_{F}\}.$ The singularity orb$(\sigma_{F})\in U_{F}$ is resolved minimally by
the proper birational map induced by the refinement $\{ \mathbb{R}_{\geq0}\,
\mathbf{u}_{j}^{(F)}+\mathbb{R}_{\geq0}\, \mathbf{u}_{j+1}^{(F)}\ \left\vert
\ 0\leq j\leq s_{F}\right.  \}$ of the fan which is composed of the cone
$\sigma_{F}$ and its faces. The exceptional divisor is $E^{(F)}:=%
{\textstyle\sum\nolimits_{j=1}^{s_{F}}}
E_{j}^{(F)},$ having%
\[
E_{j}^{(F)}:=\text{ }\overline{\text{orb}(\mathbb{R}_{\geq0}\, \mathbf{u}%
_{j}^{(F)})}\ (\cong\mathbb{P}_{\mathbb{C}}^{1}),\ \forall j\in\{1,\ldots
,s_{F}\},
\]
(i.e., the closures of the $\mathbb{T}$-orbits of the \textquotedblleft
new\textquotedblright\ rays) as its components, with self-intersection number
$(E_{j}^{(F)})^{2}=-b_{j}^{(F)}$ (see \cite[Cor. 1.18 and Prop. 1.19, pp.
23-25]{Oda}).

\begin{note}
(i) If $F\in\mathcal{F}(Q),$ and ${\boldsymbol{\eta}}_{F}\in(\mathbb{Z}%
^{2})^{\vee}$ is its unique primitive outer normal vector, we define its
\textit{local index} to be the positive integer $l_{F}:=\left\langle
{\boldsymbol{\eta}}_{F},F\right\rangle ,$ where
\[
\left\langle \cdot,\cdot\right\rangle :\text{Hom}_{\mathbb{R}}(\mathbb{R}%
^{2},\mathbb{R})\times\mathbb{R}^{2}\longrightarrow\mathbb{R}%
\]
is the usual inner product. For $F\in\mathcal{F}(Q)\mathbb{r}I_{Q}$ we have
obviously $l_{F}=1.$ For $F\in I_{Q},$ let $K(E^{(F)})$ be the \textit{local
canonical divisor }of the minimal resolution of orb$(\sigma_{F})\in U_{F}$ (in
the sense of \cite[p. 75]{Dais}). $K(E^{(F)})$ is a $\mathbb{Q}$-Cartier
divisor (a rational linear combination of $E_{j}^{(F)}$'s), and
\begin{equation}
l_{F}=\text{ min}\left\{  \left.  \xi\in\mathbb{N}\ \right\vert \ \xi
K(E^{(F)})\text{ is a Cartier divisor}\right\}  =\tfrac{q_{F}}{\text{gcd}%
(q_{F},p_{F}-1)}. \label{localind}%
\end{equation}
(ii) If $F\in I_{Q},$ denoting by $\mathfrak{m}_{X_{Q},\text{orb}(\sigma_{F}%
)}$ the maximal ideal of the local ring $\mathcal{O}_{X_{Q},\text{orb}%
(\sigma_{F})}$ of the singularity orb$(\sigma_{F}),$ and by
\[
m_{F}:=\text{dim}_{\mathbb{C}}((\mathfrak{m}_{X_{Q},\text{orb}(\sigma_{F}%
)})/(\mathfrak{m}_{X_{Q},\text{orb}(\sigma_{F})}^{2}))-1
\]
its \textit{multiplicity}, it is known (cf. \cite[Satz 2.11, p. 347]%
{Brieskorn}) that
\begin{equation}
m_{F}=2+\sum_{j=1}^{s_{F}}(b_{j}^{(F)}-2). \label{multformula}%
\end{equation}

\end{note}

\begin{lemma}
\label{MULTIND}For all $F\in I_{Q}$ we have
\[
m_{F}\leq2l_{F}.
\]
\end{lemma}

\begin{proof}
See \cite[Lemma 1.1 (iii), p. 235]{Nikulin1}.
\end{proof}

\begin{lemma}
\label{KESQUARE}For all $F\in I_{Q}$ the self-intersection number of
$K(E^{(F)})$ equals
\[
K(E^{(F)})^{2}=-\left(  \frac{2-\left(  p_{F}+\widehat{p}_{F}\right)  }{q_{F}%
}+(m_{F}-2)\right)  .
\]

\end{lemma}

\begin{proof}
Follows from \cite[Corollary 4.6, p. 96]{Dais} and formula
(\ref{multformula}).
\end{proof}

\noindent{}The minimal desingularization $\varphi:\widetilde{X}_{Q}%
\longrightarrow X_{Q}$ of $X_{Q}$ is constructed by means of the smooth
compact toric surface $\widetilde{X}_{Q}$ which is defined by the fan
\[
\widetilde{\Delta}_{Q}:=\left\{
\begin{array}
[c]{c}%
\text{ the cones }\left\{  \left.  \sigma_{F}\ \right\vert \ F\in
\mathcal{F}(Q)\mathbb{r}I_{Q}\right\}  \text{ and }\\
\left\{  \left.  \mathbb{R}_{\geq0}\, \mathbf{u}_{j}^{(F)}+\mathbb{R}_{\geq
0}\, \mathbf{u}_{j+1}^{(F)}\ \right\vert \ F\in I_{Q},\ j\in\{0,1,\ldots
,s_{F}\} \right\}  ,\\
\text{together with their faces}%
\end{array}
\right\}
\]
(refining each of the cones $\left\{  \left.  \sigma_{F}\ \right\vert \ F\in
I_{Q}\right\}  $ of $\Delta_{Q}$ as mentioned above). Furthermore, the
corresponding \textit{discrepancy divisor} equals%
\begin{equation}
K_{\widetilde{X}_{Q}}-\varphi^{\star}K_{X_{Q}}=\sum_{F\in I_{Q}}K(E^{(F)}).
\label{DISCREPANCY}%
\end{equation}
(By $K_{X_{Q}},K_{\widetilde{X}_{Q}}$ we denote the canonical divisors of
$X_{Q}$ and $\widetilde{X}_{Q},$ respectively.)

\begin{note}
By virtue of (\ref{localind}) and (\ref{DISCREPANCY}) the index $\ell$ of
$X_{Q}$ (as defined in \S \ref{INTRO}) equals%
\begin{equation}
\ell=\text{ lcm}\left\{  \left.  l_{F}\ \right\vert \ F\in\mathcal{F}%
(Q)\right\}  . \label{LCM}%
\end{equation}
(For simplicity, sometimes $\ell $ is referred as \textit{index} of $Q.$) In fact, \ if we denote by
\[
Q^{\ast}:=\left\{  \left.  \mathbf{y}\in\text{Hom}_{\mathbb{R}}(\mathbb{R}%
^{2},\mathbb{R})\ \right\vert \ \left\langle \mathbf{y},\mathbf{x}%
\right\rangle \,\leq\, 1,\ \forall\, \mathbf{x}\in Q\right\}
\]
the \textit{polar} of the polygon $Q,$ the index $\ell$ is nothing but
min$\left\{  \left.  k\in\mathbb{N}\; \right\vert \ \mathcal{V}(kQ^{\ast
})\subset\mathbb{Z}^{2}\right\}  ,$ where $kQ^{\ast}:=$ $\left\{  \left.
k\mathbf{y}\right\vert \mathbf{y}\in Q^{\ast}\right\}  .$ In other words,
$\ell$ equals the least common multiple of the (smallest) denominators of the
(rational) coordinates of the vertices of $Q^{\ast}.$
\end{note}

\section{Proof of main theorem}

\noindent{}The proof follows from suitable combination of the two upper bounds
given in Lemmas \ref{LemmaVQ} and \ref{LemmaMINDES}. (Henceforth we use freely
the notation introduced in \S \ref{PRELIM}.)

\begin{lemma}
\label{LemmaVQ}Let $X_{Q}$ be a toric log Del Pezzo surface of index $\ell
\geq1.$ Then%
\begin{equation}
\sharp(\mathcal{V}(Q))\leq4\, \text{\emph{max}}\left\{  \left.  l_{H}%
\right\vert \,H\in\mathcal{F}(Q)\right\}  +2\leq4\ell+2. \label{firstineq}%
\end{equation}
Moreover, $\sharp(\mathcal{V}(Q))= \,4\, \text{\emph{max}}\left\{  \left.
l_{H}\right\vert \,H\in\mathcal{F}(Q)\right\}  +2$, if and only if $\ell=1$,
and $Q$ is the unique hexagon \emph{(}up to lattice-equivalence\emph{)} with
one interior lattice point. This means, in particular, that for indices
$\ell\geq2$ we have%
\begin{equation}
\sharp(\mathcal{V}(Q))\leq4\ell+1. \label{nicebound}%
\end{equation}

\end{lemma}

\begin{proof}
Obviously, there exists a facet $F\in\mathcal{F}(Q)$ such that $\sum
_{\mathbf{v}\in\mathcal{V}(Q)}\mathbf{v}\in\sigma_{F}$ (this is a
\emph{special facet}, in the sense of \cite[Sect.~3]{Oebro}). In addition, since
$Q$ is two-dimensional, we have for all integers $j$:%
\[
\sharp\left\{  \left.  \mathbf{v}\in\mathcal{V}(Q)\right\vert \,\left\langle
{\boldsymbol{\eta}}_{F},\mathbf{v}\right\rangle =j\right\}  \leq2.
\]
Writing $\mathcal{V}(Q)$ as disjoint union $\mathcal{V}(Q)=\mathcal{V}_{\geq
0}^{\left(  F\right)  }(Q)\,%
{\textstyle\bigsqcup}
\,\mathcal{V}_{<0}^{\left(  F\right)  }(Q),$ where%
\[
\mathcal{V}_{\geq0}^{\left(  F\right)  }(Q):=\left\{  \left.  \mathbf{v}%
\in\mathcal{V}(Q)\right\vert \,\left\langle {\boldsymbol{\eta}}_{F}%
,\mathbf{v}\right\rangle \geq0\right\}  \text{ \ and \ \ }\mathcal{V}%
_{<0}^{\left(  F\right)  }(Q):=\left\{  \left.  \mathbf{v}\in\mathcal{V}%
(Q)\right\vert \,\left\langle {\boldsymbol{\eta}}_{F},\mathbf{v}\right\rangle
<0\right\}  ,
\]
we observe that
\[
\sharp(\mathcal{V}_{\geq0}^{\left(  F\right)  }(Q))\leq2\left(  l_{F}%
+1\right)  ,
\]
because $\left\langle {\boldsymbol{\eta}}_{F},\mathbf{v}\right\rangle
\in\{0,1,\ldots,l_{F}\}$ for all $\mathbf{v}\in\mathcal{V}_{\geq0}^{\left(
F\right)  }(Q)$. On the other hand,%
\begin{align*}
0  &  \leq\left\langle {\boldsymbol{\eta}}_{F},\sum\nolimits_{\mathbf{v}%
\in\mathcal{V}(Q)}\mathbf{v}\right\rangle =\sum\nolimits_{\mathbf{v}%
\in\mathcal{V}_{\geq0}^{\left(  F\right)  }(Q)}\left\langle {\boldsymbol{\eta
}}_{F},\mathbf{v}\right\rangle +\sum\nolimits_{\mathbf{v}\in\mathcal{V}%
_{<0}^{\left(  F\right)  }(Q)}\left\langle {\boldsymbol{\eta}}_{F}%
,\mathbf{v}\right\rangle \\
&  =\sum_{j=0}^{l_{F}}\ \sum\limits_{\left.  \{\mathbf{v}\in\mathcal{V}%
_{\geq0}^{\left(  F\right)  }(Q)\right\vert \,\left\langle {\boldsymbol{\eta}%
}_{F},\mathbf{v}\right\rangle =j\}}\left\langle {\boldsymbol{\eta}}%
_{F},\mathbf{v}\right\rangle +\sum\limits_{\mathbf{v}\in\mathcal{V}%
_{<0}^{\left(  F\right)  }(Q)}\left\langle {\boldsymbol{\eta}}_{F}%
,\mathbf{v}\right\rangle \\
&  \leq\sum_{j=0}^{l_{F}}2j+\sum\limits_{\mathbf{v}\in\mathcal{V}%
_{<0}^{\left(  F\right)  }(Q)}\left\langle {\boldsymbol{\eta}}_{F}%
,\mathbf{v}\right\rangle .
\end{align*}
This implies%
\[
a:=-\sum\nolimits_{\mathbf{v}\in\mathcal{V}_{<0}^{\left(  F\right)  }%
(Q)}\left\langle {\boldsymbol{\eta}}_{F},\mathbf{v}\right\rangle \leq
2\binom{l_{F}+1}{2}.
\]
Setting $\mu:=\sharp(\mathcal{V}_{<0}^{\left(  F\right)  }(Q))$ we examine two
cases: (i) If $\mu=2\lambda,$ for a $\lambda\in\mathbb{N},$ then
\[
\sum_{j=0}^{\lambda}2j\leq a\Longrightarrow2\binom{\lambda+1}{2}\leq
2\binom{l_{F}+1}{2}\Longrightarrow\lambda\leq l_{F}\text{ and }\mu\leq2l_{F}.
\]
(ii) If $\mu=2\lambda+1,$ for a $\lambda\in\mathbb{Z}_{\geq0},$ then
$\sum_{j=0}^{\lambda}2j+\left(  \lambda+1\right)  \leq a,$ i.e.,%
\[
2\binom{\lambda+1}{2}+\left(  \lambda+1\right)  \leq2\binom{l_{F}+1}%
{2}\Longrightarrow\lambda\leq l_{F}-1\text{ and }\mu\leq2l_{F}-1.
\]
Hence,%
\[
\sharp(\mathcal{V}(Q))=\sharp(\mathcal{V}_{\geq0}^{\left(  F\right)
}(Q))+\sharp(\mathcal{V}_{<0}^{\left(  F\right)  }(Q))\leq2\left(
l_{F}+1\right)  +\mu
\]%
\begin{equation}
\leq2\left(  l_{F}+1\right)  +2l_{F}=4l_{F}+2\leq4\,\text{max}\left\{  \left.
l_{H}\right\vert \,H\in\mathcal{F}(Q)\right\}  +2, \label{Approx}%
\end{equation}
with the latter upper bound $\leq4\ell+2$ (by (\ref{LCM})), giving the
inequality (\ref{firstineq}). Finally, we deal with the case of equality:
Suppose that $\sharp(\mathcal{V}(Q))=4\ell^{\prime}+2$, where
\[
\ell^{\prime}:=\max\left\{  \left.  l_{H}\right\vert \,H\in\mathcal{F}%
(Q)\right\}  .
\]
From (\ref{Approx}) we see that $\mu=2l_{F}$, and $\lambda=l_{F}=\ell^{\prime
}$. Therefore, by the equalities in (i) we have for the integers
$j=-\ell^{\prime},\ldots,0,\ldots,\ell^{\prime}$:
\begin{equation}
\sharp\left\{  \left.  \mathbf{v}\in\mathcal{V}(Q)\right\vert \,\left\langle
{\boldsymbol{\eta}}_{F},\mathbf{v}\right\rangle =j\right\}  =2.
\label{gleichheit}%
\end{equation}
In particular, $0=\left\langle {\boldsymbol{\eta}}_{F},\sum
\nolimits_{\mathbf{v}\in\mathcal{V}(Q)}\mathbf{v}\right\rangle $, i.e.,
$\sum_{v\in\mathcal{V}(Q)}\mathbf{v}=\mathbf{0}$. Hence, the previous argument
holds for \textit{any} facet. Now let $F^{\prime}$ be another facet of $Q$
having a common vertex, say $\mathbf{v,}$ with $F.$ If $\mathcal{V}%
(F)=\{\mathbf{u},\mathbf{v}\}$ and $\mathcal{V}(F^{\prime})=\{\mathbf{v}%
,\mathbf{w}\},$ then applying (\ref{gleichheit}) for \textit{both} $F$ and
$F^{\prime}$ we get $\left\langle {\boldsymbol{\eta}}_{F},\mathbf{w}%
\right\rangle =\ell^{\prime}-1$ and $\left\langle {\boldsymbol{\eta}%
}_{F^{\prime}},\mathbf{u}\right\rangle =\ell^{\prime}-1$. This implies
$\ell^{\prime}=1=\ell$, since otherwise the primitive vertex $\mathbf{v}$
equals $(\ell^{\prime}/(\ell^{\prime}-1))(\mathbf{w}+\mathbf{u}-\mathbf{v})$,
a contradiction. Consequently, $Q$ has to be the unique hexagon (up to
lattice-equivalence) with just one interior lattice point (see \cite[Proposition 2.1]{Nill}).
\end{proof}

\begin{lemma}
\label{LemmaMINDES}If $X_{Q}$ is a toric log Del Pezzo surface of index
$\ell\geq2$ and $\widetilde{X}_{Q}\overset{\varphi}{\longrightarrow}X_{Q}$ its
minimal desingularization, then%
\begin{equation}
\rho(\widetilde{X}_{Q})<2\, \sharp(I_{Q}\mathbb{r}\breve{I}_{Q})(\ell-1)
-\frac{1}{\ell} \; \sharp(\mathcal{V}(Q)) +10. \label{secineq}%
\end{equation}

\end{lemma}

\begin{proof}
By Noether's formula and (\ref{DISCREPANCY}) we deduce
\[
\rho(\widetilde{X}_{Q})=10-K_{\widetilde{X}_{Q}}^{2}=10-K_{X_{Q}}^{2}%
-\sum_{F\in I_{Q}}K(E^{(F)})^{2}.
\]
Since $-\ell K_{X_{Q}}$ is an ample Cartier divisor on $X_{Q},$ we can compute
by \cite[Proposition 2.10, p. 79]{Oda} its self-intersection number:
\[
(-\ell K_{X_{Q}})^{2}=2\text{\thinspace area}(\ell Q^{\ast})\Longrightarrow
K_{X_{Q}}^{2}=\frac{2}{\ell^{2}}\,\text{area}(\ell Q^{\ast})=2\text{\thinspace
area}(Q^{\ast}).
\]
For any facet $H$ of $\ell Q^{\ast}$ the primitive outer normal vector is given by
some vertex of $Q$, i.e., the lattice distance of $H$ from $\mathbf{0}$ equals
$\ell$. This implies
\[
\text{area}(\ell Q^{\ast}) \geq\frac{1}{2} \ell\;\sharp(\mathcal{F}(\ell
Q^{\ast})) = \frac{1}{2} \ell\;\sharp(\mathcal{V}(Q)).
\]
Hence,
\[
-K_{X_{Q}}^{2} = -\frac{2}{\ell^{2}} \;\text{area}(\ell Q^{\ast}) \leq
-\frac{1}{\ell} \;\sharp(\mathcal{V}(Q)).
\]
On the other hand, by Lemma \ref{KESQUARE} we infer that%
\[
-\sum_{F\in I_{Q}}K(E^{(F)})^{2}=\sum_{F\in I_{Q}}\left(  \frac{2-\left(
p_{F}+\widehat{p}_{F}\right)  }{q_{F}}+(m_{F}-2)\right)  .
\]
Taking into account that $m_{F}=2$ for all $F\in\breve{I}_{Q},$ and that
$p_{F}+\widehat{p}_{F}\geq2$ for all $F\in I_{Q},$ which is valid as equality
only for $p_{F}=\widehat{p}_{F}=1,$ i.e., whenever $F\in\breve{I}_{Q},$ we
obtain%
\[
-\sum_{F\in I_{Q}}K(E^{(F)})^{2}=-\sum_{F\in I_{Q}\mathbb{r}\breve{I}_{Q}%
}K(E^{(F)})^{2}<\sum_{F\in I_{Q}\mathbb{r}\breve{I}_{Q}}(m_{F}-2)
\]%
\begin{align*}
&  \leq\sharp(I_{Q}\mathbb{r}\breve{I}_{Q})\,\,\text{max}\left\{  \left.
m_{F}-2 \;\right\vert \; F\in I_{Q}\mathbb{r}\breve{I}_{Q}\right\} \\
&  \leq\,\sharp(I_{Q}\mathbb{r}\breve{I}_{Q})\,\,\text{max}\left\{  \left.
2(l_{F}-1) \;\right\vert \; F\in I_{Q}\mathbb{r}\breve{I}_{Q}\right\}
\leq2\,\sharp(I_{Q}\mathbb{r}\breve{I}_{Q})(\ell-1),
\end{align*}
where the last but one inequality follows from Lemma \ref{MULTIND}. Thus,
$\rho(\widetilde{X}_{Q})$ is strictly smaller than the sum $10 -\sharp
(\mathcal{V}(Q))/\ell+ 2\,\sharp(I_{Q}\mathbb{r}\breve{I}_{Q})(\ell-1).$
\end{proof}

\noindent{}\textit{Proof of Theorem \ref{main}}. If $\ell=1,$ then
$\rho(\widetilde{X}_{Q})\leq7$ by the known classification of the reflexive
polygons (see \cite{KS} or \cite[Proposition 2.1]{Nill}). If $\ell\geq2,$
applying (\ref{nicebound}) and (\ref{secineq}), and the inequality
$\sharp(I_{Q}\mathbb{r}\breve{I}_{Q})\leq\sharp(\mathcal{V}(Q)),$ we get
\[
\rho(\widetilde{X}_{Q})<2\,\sharp(I_{Q}\mathbb{r}\breve{I}_{Q})(\ell
-1)-\frac{1}{\ell}\;\sharp(\mathcal{V}(Q))+10
\]%
\[
\leq\sharp(\mathcal{V}(Q))\left(  2(\ell-1)-\frac{1}{\ell}\right)
+10\leq\left(  4\ell+1\right)  \left(  2(\ell-1)-\frac{1}{\ell}\right)  +10,
\]
i.e., $\rho(\widetilde{X}_{Q})<8\ell^{2}-6\ell+4-\frac{1}{\ell}$, which yields the bound for $\ell \geq 2$.
\hfill{}$\square$

\section{Discussion, improvements and examples}

\noindent{} First, let us note that from the proof of Theorem \ref{main} we derive a \emph{linear} upper bound on
$\rho(\widetilde{X}_{Q})$, if the number of vertices of $Q$ is \emph{fixed}. It is
therefore natural to ask for an example of an infinite family $\{Q_{i}\}$ of
LDP-polygons with increasing number of vertices, for which $\rho(\widetilde
{X}_{Q_{i}})$ exhibits a non-linear growth with respect to the indices of its members.
To the best knowledge of the authors, this seems to be an open question.

Now, in some specific cases we can further improve the bound (\ref{mainbound}).
If $Q$ is an LDP-polygon and $F\in I_{Q},$ then, according to (\ref{localind}%
), there is a positive integer $\beta_{F}$ such that%
\[
p_{F}-1=\beta_{F}\cdot\frac{q_{F}}{l_{F}}\Longrightarrow l_{F}\,(p_{F}%
-1)=\beta_{F}\,q_{F}.
\]
Since $l_{F}(p_{F}-1)<l_{F}(q_{F}-1)<l_{F}q_{F},$ we have $\beta_{F}%
\in\{1,\ldots,l_{F}-1\}.$ In Proposition~\ref{SPECIALCASE} we construct a
better upper bound for $\rho(\widetilde{X}_{Q})$ provided that $\beta_{F}$
takes one of the extreme values $1,l_{F}-1,$ and $l_{F}^{2}\mid q_{F}$ for all
$F\in I_{Q}\mathbb{r}\breve{I}_{Q}.$

\begin{proposition}
\label{SPECIALCASE}Let $Q$ be an LDP-polygon such that $X_{Q}$ has index
$\ell\geq2.$ Suppose that for all $F\in I_{Q}\mathbb{r}\breve{I}_{Q}$ the
following conditions are satisfied\emph{:\smallskip} \newline\emph{(i)}
$\beta_{F}\in\{1,l_{F}-1\},$ and\smallskip\ \newline\emph{(ii)} $l_{F}^{2}\mid
q_{F}.$ \noindent Then%
\begin{equation}
\rho(\widetilde{X}_{Q})\leq4\ell^{2}-3\ell+4. \label{LINUB}%
\end{equation}

\end{proposition}

\begin{proof}
For $F\in I_{Q}\mathbb{r}\breve{I}_{Q}$ define $\xi_{F}:=\frac{q_{F}}{\ell
^{2}}.$ If $\beta_{F}=1,$ then $\frac{q_{F}}{q_{F}-p_{F}}$ equals
\[
\left\{
\begin{array}
[c]{ll}%
1+\frac{1}{(l_{F}-2)+\frac{1}{l_{F}+1}}=[\![\underset{\left(  l_{F}-2\right)
\text{\emph{-}times}}{\underbrace{2,...,2}},l_{F}+2]\!], & \text{if }\xi
_{F}=1,\\
\, & \\
1+\frac{1}{l_{F}-2+\frac{1}{1+\frac{1}{\xi_{F}-1+\frac{1}{l_{F}}}}%
}=[\![\underset{\left(  l_{F}-2\right)  \text{{\small -times}}}{\underbrace
{2,..,2}},3,\underset{\left(  \xi_{F}-2\right)  \text{{\small -times}}%
}{\underbrace{2,..,2}},\ell+1]\!], & \text{if }\xi_{F}\geq2,
\end{array}
\right.
\]
(cf. \cite[Proposition 3.1, pp. 83-84]{Dais}), $\widehat{p}_{F}=q_{F}-l_{F}%
\xi_{F}+1,$ and
\[
m_{F}-2=%
{\textstyle\sum\limits_{j=1}^{s_{F}}}
(b_{j}^{(F)}-2)=l_{F},\ \forall F\in I_{Q}\mathbb{r}\breve{I}_{Q}.
\]
Correspondingly, if $\beta_{F}=l_{F}-1,$ then $\frac{q_{F}}{q_{F}-p_{F}}$
equals%
\[
\left\{
\begin{array}
[c]{ll}%
(l_{F}+1)+\frac{1}{l_{F}-1}=[\![l_{F}+2,\underset{\left(  l_{F}-2\right)
\text{\emph{-}times}}{\underbrace{2,...,2}}]\!], & \text{if }\xi_{F}=1,\\
\, & \\
l_{F}+\frac{1}{(\xi_{F}-1)+\frac{1}{1+\frac{1}{l_{F}-1}}}=[\![\ell
+1,\underset{\left(  \xi_{F}-2\right)  \text{{\small -times}}}{\underbrace
{2,...,2}},3,\underset{\left(  l_{F}-2\right)  \text{{\small -times}}%
}{\underbrace{2,...,2}}]\!], & \text{if }\xi_{F}\geq2,
\end{array}
\right.
\]
$\widehat{p}_{F}=l_{F}\xi_{F}+1,$ and $m_{F}-2=%
{\textstyle\sum\limits_{j=1}^{s_{F}}}
(b_{j}^{(F)}-2)=l_{F},\ \forall F\in I_{Q}\mathbb{r}\breve{I}_{Q}.$ Thus,%
\begin{align*}
-\sum_{F\in I_{Q}\mathbb{r}\breve{I}_{Q}}K(E^{(F)})^{2}  &  =\sum_{F\in
I_{Q}\mathbb{r}\breve{I}_{Q}}\left(  \tfrac{2-\left(  p_{F}+\widehat{p}%
_{F}\right)  }{q_{F}}+(m_{F}-2)\right) \\
&  =\sum_{F\in I_{Q}\mathbb{r}\breve{I}_{Q}}(l_{F}-1)\leq\sharp(I_{Q}%
\mathbb{r}\breve{I}_{Q})(\ell-1).
\end{align*}
Since $\sharp(I_{Q}\mathbb{r}\breve{I}_{Q})\leq\sharp(\mathcal{V}(Q)),$
applying Lemma \ref{LemmaVQ} and the reasoning used in the proof of Lemma
\ref{LemmaMINDES}, we get%
\begin{align*}
\rho(\widetilde{X}_{Q})  &  <\sharp(I_{Q}\mathbb{r}\breve{I}_{Q}%
)(\ell-1)-\frac{1}{\ell}\;\sharp(\mathcal{V}(Q))+10\\
&  \leq\sharp(\mathcal{V}(Q))\left(  \ell-1-\frac{1}{\ell}\right)
+10\leq(4\ell+1)\left(  \ell-1-\frac{1}{l}\right)  +10.
\end{align*}
The upper bound (\ref{LINUB}) follows from this inequality.
\end{proof}

By \cite[Lemma 6.9, p. 107]{Dais} we see that the conditions (i), (ii) in
Proposition \ref{SPECIALCASE} are automatically satisfied for all toric log
Del Pezzo surfaces of index $\ell=2.$ Hence, the upper bound $14$ improves
noticeably (\ref{mainbound}) (which equals $23$ in this case). In fact, for
$\ell=2,$ it can be shown (though, at the cost of passing through ad hoc
classification results for the corresponding LDP-polygons) that the
\textit{sharp} upper bound equals $10$.

Finally, in Proposition \ref{LDPRANK1} we classify those LDP-triangles of {\em arbitrary} index, whose
toric log Del Pezzo surfaces have at most one singularity. Somehow surprisingly, the Picard number
of their minimal desingularizations is bounded; moreover, it takes always the {\em smallest} possible value, namely $2$. Note that from Lemma \ref{LemmaMINDES}
one only derives that the Picard numbers behave {\em at most linearly} with respect to the index, once the number of non-Gorenstein singularities
$\sharp(I_{Q}\mathbb{r}\breve{I}_{Q})$ is fixed.

\begin{lemma}
\label{LEMMAQP}If $X_{Q}$ is a toric log Del Pezzo surface with Picard number
$\rho(X_{Q})=1$ \emph{(}i.e, if $Q$ is an LDP-\emph{triangle) }and
$\sharp(I_{Q})=1,$ then $Q$ is lattice-equivalent to the triangle $Q_{p}$
having $\binom{1}{0},\binom{p}{p+1}$ and $\binom{-1}{-1}$ as its vertices, for
some positive integer $p.$
\end{lemma}

\begin{proof}
If $I_{Q}=\{F\},$ setting $p:=p_{F}$ and $q:=q_{F},$ there is a unimodular
transformation mapping $\mathbf{n}^{(F)}$ onto $\mathbf{n}_{1}:=\binom{1}{0},$
$\mathbf{n}^{\prime(F)}$ onto $\mathbf{n}_{2}:=\binom{p}{q},$ and the third
vertex of $Q$ onto an $\mathbf{n}_{3}=\binom{x_{1}}{x_{2}}$ which belongs
necessarily to the set $\left\{  \binom{x_{1}}{x_{2}}\in\mathbb{Z}%
^{2}\ \left\vert \ \frac{q}{p}x_{1}<x_{2}<0\right.  \right\}  .$ Since
$\left\vert \det(\mathbf{n}_{2},\mathbf{n}_{3})\right\vert =$ $\left\vert
\det(\mathbf{n}_{3},\mathbf{n}_{1})\right\vert =1,$ we have $x_{2}=-1$ and
$x_{1}=-\frac{p+1}{q}.$ Hence, $q\mid p+1,$ which implies $q=p+1\, \ $(because
$p<q$).
\end{proof}

\begin{proposition}
\label{LDPRANK1}Let $X_{Q}$ be a toric log Del Pezzo surface which has Picard
number $\rho(X_{Q})=1,$ arbitrary index $\ell\geq1,$ and $\sharp(I_{Q})=1.$
For $\ell$ odd $\geq3$ we have either $X_{Q}\cong X_{Q_{\ell-1}}$ or
$X_{Q}\cong X_{Q_{2\ell-1}},$ whereas for $\ell\in\{1\} \cup2\mathbb{Z}$ we
have $X_{Q}\cong X_{Q_{2\ell-1}}.$ Furthermore, for all $\ell\geq1,$ the
Picard number of the rational surface $\widetilde{X}_{Q}$ obtained by the
minimal resolution of the singularity of $X_{Q}$ equals
\[
\rho(\widetilde{X}_{Q})=2.
\]

\end{proposition}

\begin{proof}
By Lemma \ref{LEMMAQP}, the LDP-triangle\textit{ }$Q$ is lattice-equivalent to
$Q_{p},$ for some positive integer $p.$ Since $q=p+1$ and gcd$(p+1,p-1)\in
\{1,2\}$, the index $\ell$ of $X_{Q}\cong X_{Q_{p}}$ equals $\frac{p+1}{2}$
whenever $p$ is odd and $p+1$ whenever $p$ is even (see (\ref{LCM})). This
bears out our first assertion. On the other hand, since $\widetilde{\Delta
}_{Q_{p}}$ is obtained from $\Delta_{Q_{p}}$ by adding just one new ray
(namely $\mathbb{R}_{\geq0}\binom{1}{1}$), we have%
\[
\rho(\widetilde{X}_{Q})=\rho(\widetilde{X}_{Q_{p}})=\sharp\{ \text{rays of
\ }\widetilde{\Delta}_{Q_{p}}\}-2=4-2=2,
\]
(cf. \cite[Corollary 2.5, p. 74]{Oda}). Thus, the second assertion is also true.
\end{proof}

It would be interesting to generalize this result by regarding LDP-{\em polygons} of arbitrary index, whose
toric log Del Pezzo surfaces have at most one singularity.

\end{document}